\documentclass[notitlepage,draft,10pt,times]{article}

\usepackage{amssymb}
\usepackage{amsmath,amsthm,amssymb, mathrsfs,enumerate,wasysym}
 \usepackage{xypic}
 \usepackage[left=1.5cm,top=2cm,right=1.5cm]{geometry}
\usepackage[all]{xy}
\usepackage[OT2,T1]{fontenc}
\usepackage{textcomp}
\usepackage{fancyhdr}
\DeclareSymbolFont{cyrletters}{OT2}{wncyr}{m}{n}
\DeclareMathSymbol{\Sha}{\mathalpha}{cyrletters}{"58}

\newcommand{\Q}{\mathbb{Q}}

\newcommand{\Z}{\mathbb{Z}}

\newcommand{\links}{\left(\begin{array}{cc}}
\newcommand{\rechts}{\end{array}\right)}
\newcommand{\bai}{\left[\begin{array}{cc}}
\newcommand{\dai}{\end{array}\right]}
\newcommand{\hidari}{\left(\begin{array}{c}}
\newcommand{\migi}{\end{array}\right)}

\newcommand{\Gal}{{\mathrm{Gal}}}

\renewcommand{\phi}{{\varphi}}

\newtheorem{theorem}{Theorem}[section]

\newtheorem{open problem}[theorem]{Open Problem}

\theoremstyle{definition}

\begin{document}



\title{Elliptic Curves with Supersingular Reduction over $\Gamma$-extensions}
\date{}
\author{A.G. Nasybullin}
\maketitle

\thispagestyle{empty}
\renewcommand{\headrulewidth}{0.0pt}
\thispagestyle{fancy}
\fancyhead[c]{In the Moscow Mathematical Society}







\begin{section}{}
Let $p$ be a prime number, $k_0$ a finite extension of the rationals $\Q$, $k_\infty/k_0$ a Galois extension with [Galois] group $\Gamma$ isomorphic to the group of $p$-adic integers $\Z_p$. Put $\Gamma_n:=\Gamma^{p^n}, k_n:=k_\infty^{\Gamma_n}$. Let $E$ be an elliptic curve over $\Q$ with supersingular reduction at $p$, $E(k_n)$ the $k_n$-rational points of $E$, and $\Sha_n^{(p)}$ the $p$-component of the Shafarevich-Tate group of the curve $E\otimes k_n$.
\begin{theorem}\label{firsttheorem}
We assume the following conditions:
\begin{enumerate}[(a)]
\item $p $ is not $2$ and does not divide the number of the rational connected components of bad reduction of the curve $E\otimes k_0$.
\item For all places $v$ of $k_0$ dividing $p$, the completion $k_{0,v}$ is unramified over the field of the $p$-adic numbers $\Q_p$, and its degree over $\Q_p$ is not divisible by $4$.
\item The $\Gamma$-extension $k_\infty/k_0$ is cyclotomic, i.e.
\[\displaystyle k_\infty \subset \bigcup_{n=1}^\infty k_0(\sqrt[p^n]{1})\]
\end{enumerate}
Then, if $E(k_0)$ is finite and $\Sha_0^{(p)}=0$, the groups $E(k_n),E(k_\infty)$ and $\Sha_n^{(p)}$ are finite and
\[\log_p [\Sha_n^{(p)}]=[k_0:\Q]\left(\left[\frac{p^{n+1}}{p^2-1}\right]-\left[\frac{n+1}{2}\right]\right).\]
\end{theorem}

We denote by $a_p$ the trace of the Frobenius automorphism of the reduction of $E \mod p$. Note that $E\mod p$ is supersingular if and only if it is non-singular and $p$ divides $a_p$. Consequently, $a_p=0$ for $p>3$ and $a_p=0, \pm p$ when $p=2,3$.

\begin{theorem}\label{middletheorem}Suppose that $k_0/\Q$ is abelian and $k_\infty/k_0$ is cyclotomic. Then:
\begin{enumerate}[(a)]
\item There are integers $\rho^{(0)}, \rho^{(1)}\geq 0$, equal for $a_p\neq 0$,  such that for all sufficiently large $ n\equiv s \mod 2$ $ (s=0,1)$,
\[\mathrm{rk \, } E(k_n)+\mathrm{cork}\, \Sha_n^{(p)}-\mathrm{rk}\, E(k_{n-1})-\mathrm{cork}\,\Sha_{n-1}^{(p)}=\rho^{(s)}(p^n-p^{n-1}),\]
where $\mathrm{rk \, } E(k_n)$ is the rank of $E(k_n)$ and $\mathrm{cork}\, \Sha_n^{(p)}$ is the corank of $\Sha_n^{(p)}$;
\item if $a_p\neq0$ and the degree $[k_0:\Q]$ divides a number of the form $(p^l+1)p^m$, then $\mathrm{rk} \, E(k_n)$ stabilizes, and consequently $E(k_\infty)$ is finitely generated;
\item if $E(k_0)$ and $\Sha_0^{(p)}$ are finite, and for $a_p=0$ we have the condition (b) of Theorem \ref{firsttheorem}, then $\rho^{(0)}=\rho^{(1)}=0$, i.e. $\mathrm{rk \, } E(k_n)$ and $\mathrm{cork}\, \Sha_n^{(p)}$ stabilize;
\item if $\rho^{(0)}=\rho^{(1)}=0$, then there are integers $\mu^{(s)},\delta^{(s)}\geq 0, \lambda^{(s)}(s=0,1)$ such that $\delta^{(0)}=\delta^{(1)}=0$ for $a_p=0$ and for all sufficiently large $n\equiv s \mod 2$,
\[\log_p[\Sha_n^{(p)}]-\log_p[\Sha_{n-1}^{(p)}]=\mu^{(s)}(p^n-p^{n-1})+([k_0:\Q]-\delta^{(s)})\left[{p^n}\over{p+1}\right]+\delta^{(s)}\left[{p^{n-1}}\over{p+1}\right]+\lambda^{(s)},\]
where $\Sha_n^{(p)}$ is the cotorsion of $\Sha_n^{(p)}$.
\end{enumerate}
\end{theorem}
Denote by $T_n$ the set of places of the field $k_n$ dividing $p$ and ramified in $k_\infty$.
\begin{theorem}\label{lasttheorem}Suppose that $a_p=0$ and that for all $n$ and $v\in T_n$, the extensions $k_{n,v}/\Q_p$ are abelian. Then there are integers $\rho^{(s)}, r^{(s)}, \nu^{(s)}, \mu^{(s)},\lambda^{(s)}\,(s=0,1), \,\mu_i \,(i=1,2,\cdots),$ satisfying the relations
\[\rho^{(s)}\geq r^{(s)}\geq\nu^{(s)}\geq0,\qquad \rho^{(0)}-r^{(0)}=\rho^{(1)}-r^{(1)},\]
\[\mu_1\geq\mu_2\geq\cdots\geq0,\mu_i=0\text{ for } i>\min(\nu^{(0)},\nu^{(1)}),\]
\[\mu^{(s)}\geq0, \, r^{(0)}+r^{(1)}\leq r, \text{ where } r=\displaystyle\sum_{v\in T_0}[k_{0,v}:\Q_p]\leq[k_0:\Q],\]
and such that for sufficiently large $n\equiv s\mod 2$ the following assertions hold:
\begin{enumerate}[(a)]
\item $\mathrm{rk \, } E(k_n)+\mathrm{cork}\, \Sha_n^{(p)}-\mathrm{rk}\, E(k_{n-1})-\mathrm{cork}\,\Sha_{n-1}^{(p)}=\rho^{(s)}(p^n-p^{n-1});$
\item 
$\log_p[\Sha_n^{(p)}]-\log_p[\Sha_{n-1}^{(p)}]=$
\[\mu^{(s)}(p^n-p^{n-1})+(r-r^{(s)}+\nu^{(s)})\left[{p^n}\over{p+1}\right]-\displaystyle \sum_{i=1}^{\nu^{(s)}}\left[{p^{n-\mu_i}}\over{p+1}\right]-\sum_{i=1}^{r^{(1-s)}}\left[{p^{n-\mu_i}}\over{p+1}\right]+\lambda^{(s)};\]
\item if $r^{(s)}=0$, then $\mathrm{rk} \, E(k_n)=\text{rk} \, E(k_{n-1});$
\item if $\mathrm{cork}\,\Sha_n^{(p)}$ stabilizes, then $\mathrm{rk}\,B_n-\mathrm{rk}\,B_{n-1}=r^{(s)}(p^n-p^{n-1}),$ where $B_n$ is the image of $E(k_n)\otimes \Z_p\rightarrow \displaystyle \sum_{v\in T_n} E(k_{n,v})^{(p)}$ and $\mathrm{rk}\, B_n$ is the rank of $B_n$ over $\Z_p$.
\end{enumerate}
\end{theorem}

In the case of nonsupersingular reduction, the behavior of the groups $E(k_n)$ and $\Sha_n^{(p)}$ has been investigated by B. Mazur (see \cite{manin}, \cite{mazur}). One of the main points of his research is the description of the $\Gamma$-modules $E(k_{n,v})^{(p)}$ for $v\in T_n$. Analogously, the proofs of Theorems \ref{firsttheorem}, \ref{middletheorem}, and \ref{lasttheorem} are based on the theorem in the following paragraph.
$\,$
\end{section}

\thispagestyle{empty}
\renewcommand{\headrulewidth}{0.0pt}
\thispagestyle{fancy}
\fancyhead[c]{In the Moscow Mathematical Society}

\begin{section}{The Local Group of Points}
Let $E$ be an elliptic curve over $\Q_p$, $E \mod p$ be supersingular, $a_p$ the trace of the Frobenius automorphism on the reduction $E \mod p$. For any abelian extension $K/\Q_p$, set $K_n:=K\cap \Q_p^{nr}(\zeta_n)$, where $n=-1,0,1,\cdots;\Q_p^{nr}$ denotes the maximal unramified extension of $\Q_p$, and $\zeta_n$ a primitive root of unity of degree $p^{n+1}$ if $p\neq2$, and of degree $p^{n+2}$ if $p=2$. We will denote by $m(K)$ the smallest $n$ for which $K_n=K$.
\begin{theorem}Let $K/\Q_p$ be  a finite abelian extension with [Galois] group $G=\Gal(K/\Q_p)$. Then the $\Z_p[G]$-module $E(K)^{(p)}$ is free of $p$-torsion and has a system of generators $\{e_n|n=-1,-,\cdots,m(K)\}$, all of whose relations can be derived from the following:
\[e_n\in E(K_n),\]
\[\mathrm{Nor}_{n/{n-1}}e_n=a_pe_{n-1}-e_{n-2}\;\;\;\;\;\;\;(n\geq 2),\]
\[\mathrm{Nor}_{1/0}e_1=\begin{cases}a_pe_0-[K_0(\zeta_0):K_0]e_{-1},&(p\neq 2),\\ a_pe_0-[K_0(\zeta_0):K_0](a_p-F-F^{-1})e_{-1},&(p=2),
\end{cases}\]
\[\mathrm{Nor}_{0/-1}e_0=\begin{cases}(a_p-F-F^{-1})e_{-1},&(p\neq 2),\\ (a_p^2-a_pF-a_pF^{-1}-1)e_{-1},&(p=2),
\end{cases}\]
where $\mathrm{Nor}_{n/{n-1}}:E(K_n)\rightarrow E(K_{n-1})$ is the norm homomorphism and $F\in \Gal(K_{-1}/\Q_p)$ is the Frobenius automorphism.
\end{theorem}
The author would like to thank Yu. I. Manin for posing the problem and his constant interest in working on it, and V. G. Berkovich for helpful discussions.
\end{section}
\hfill
\textit{  (translated by Igor Minevich and Florian Sprung.) }

\begin{flushright}
Received 20 October 1976
\end{flushright}

\thispagestyle{empty}
\renewcommand{\headrulewidth}{0.0pt}
\thispagestyle{fancy}
\fancyhead[c]{In the Moscow Mathematical Society}
\end{document}